\documentstyle{amsppt}
\magnification\magstep1 \NoRunningHeads \TagsOnRight \vsize 23 true cm \hsize 15 true cm
\hcorrection{0.6 true cm}

\centerline{\bf SIMPLE BOL LOOPS}
\bigskip\bigskip
\centerline {\bf E.K. Loginov}\footnote ""{E-mail address: loginov\@ivanovo.ac.ru}
\medskip
\centerline {\sl Ivanovo State University} \centerline {\sl Ermaka St.~39, Ivanovo, 153025,
Russia}
\bigskip\bigskip

{\bf Abstract.} In this paper we investigate the Bol loops and connected with them groups. We
prove an analog of the Doro's theorem for Moufang loops and find a criterion for simplicity of
Bol loops. One of the main results obtained is the following: if the right multiplication
group of a connected finite Bol loop $S$ is a simple group, then $S$ is a Moufang loop.
\bigskip\bigskip

{\bf 1. Introduction}
\medskip
Recall that a {\it loop} is a binary system $S$ with unit element such that the equations
$ax=b$ and $ya=b$ are uniquely solvable for all $a,b\in S$. The {\it Bol loops} are
distinguished from the class of all loops by the identity
$$
((xy)z)y=x((yz)y)
$$
A Bol loop is called a {\it Moufang loop} if it satisfies the identity
$$
y(z(yx))=(y(zy))x.
$$
These identities arose in connection with the so-called Bol closure conditions in web theory.
The example $x\cdot y=x-y$ on $\Bbb R$ shows that both identities in a quasigroup are
independent and do not imply the existence of an identity element. As well known (see,
e.g.,~[1]), every Bol loop is monoassociative and right invertible, and a left invertible Bol
loop is Moufang.
\par
A Bol loop is {\it simple} if it has no (nontrivial) proper homomorphic images, or
equivalently, if it has no proper normal subloops. Let $S$ be a finite simple Bol loop. Of
course if $S$ is associative, then $S$ is a simple group. The classification of finite simple
groups is known~[2]. This is either a cyclic group of prime order or an alternating group of
degree $n\geqslant 5$ or a group of Lie type or one of the twenty-six sporadic groups. All
finite simple nonassociative Moufang loops are also known~[3,4]. Every such loop is a
homomorphic image of the loop of all invertible elements of split Cayley-Dickson algebra over
$GF(q)$.
\par
It is generally agreed that the most significant open problem in loop theory today is the
existence of finite, simple Bol loops which are not Moufang. While we certainly have not
solved these problems, we take a step in this direction. Let $S$ be a Bol loop and $Gr_{r}(S)$
be its right multiplication group. We'll call the loop $S$ {\it strongly simple} if the group
$Gr_{r}(S)$ is simple. It is easy to prove~(see~[5]) that any strongly simple loop is simple.
If $S$ is infinite strongly simple Bol loop, then it may be non-Moufang. The corresponding
example was constructed in~[6]. If $S$ is finite, then we have the following assertion:
\medskip
{\bf Theorem.} {\sl Every  connected finite strongly simple Bol loop is Moufang. In
particular, if the right multiplication group $Gr_{r}(S)$ of a Bol loop $S$ is a finite simple
group of Lie type, then $S$ is a finite simple Moufang loop.}
\medskip
The present paper is organized as follows. In Sec. 2 we consider the typical situation of the
abstract theory of symmetric spaces: a group $G$ equipped with an involutory automorphism and
the structure of the symmetric space on $G$ defined by this pair. To this symmetric space one
naturally assigns a symmetric space $S$, and sufficient conditions are found under which the
natural binary composition transforms $S$ into a Bol loop (or a Moufang loop). In Sec. 3,
among the groups that are lead to a given Bol loop by using this construction, we choose an
universal group from which the other groups are constructed as its homomorphic images. This
makes it possible to give a criterion for simplicity of Bol loops $S$. In the last section, we
prove the theorem.
\bigskip

{\bf 2. A construction of Bol loops}
\medskip

Recall~(see, e.g.,~[7]) that a binary system $S$ is called a {\it symmetric space} if it
satisfies the identities
$$
\gather x.x=x,\\
x.(x.y)=y,\\
x.(y.z)=(x.y).(x.z)
\endgather
$$
A symmetric space with the distinguished point $e\in S$ is called a {\it punctured symmetric
space}, and the point $e$ is referred to as the basic point of the space $S$.
\par
The isomorphisms and homomorphisms of symmetric spaces are defined in the usual way; in the
case of punctured spaces, it is assumed in  addition that any homomorphism takes the basic
point to the basic point. By induction one can define the powers of an element $x$ (with
respect to the basic point $e$): $x^{0}=e$, $x^{1}=x$, $x^{n+2}=x.(e.x^{n})$, and
$x^{-n}=e.x^{n}$. It is obvious that the map $x\to x^{-1}$ is an automorphism of a symmetric
space.
\par
Suppose $G$ is a group admitting an involutory automorphism $\sigma$, $I_{\sigma}$ is the set
of fixed points of $\sigma$, and $G_{\sigma}=\{x^{-1}x^{\sigma}\mid x\in G\}$. It is easy to
prove, that the set $G_{\sigma}$ and the set of right cosets $G/I_{\sigma}$ of the subgroup
$I_{\sigma}$ in the group $G$ is a symmetric spaces with respect to the products
$$
\align
x.y&=xy^{-1}x,\\
I_{\sigma}x\cdot I_{\sigma}y&=I_{\sigma}y^{\sigma}(x^{\sigma})^{-1}x.
\endalign
$$
Moreover, the mapping $\psi: G/I_{\sigma}\to G_{\sigma}$ given by the relation
$\psi(I_{\sigma}x)=x^{-1}x^{\sigma}$ defines an isomorphism $G/I_{\sigma}\simeq G_{\sigma}$ of
the symmetric spaces.
\par
Let $N$ be the subgroup of $G$ generated by the set $G_{\sigma}$. Suppose that there is a
group homomorphism $\varphi:N\to G$ from $N$ into $G$ such that \roster
\item"{$1)$}"
$G=\langle x, x^{\sigma}\mid x\in\varphi (G_{\sigma})\rangle$,
\item"{$2)$}"
$\varphi (\varphi(x)^{-1}\varphi(x)^{\sigma})=\varphi(x)$ for any $x\in G_{\sigma}$,
\item"{$3)$}"
every coset of $I=I_{\sigma}\cap\varphi(N)$ in the group $\varphi(N)$ contains at least one
element of $\varphi(G_{\sigma})$.
\endroster
(Here and everywhere below, the notation $G=\langle M\rangle$ means that the group $G$ is
generated by the set $M$). In this case, we say that $G$ is a {\it group with the
$(\sigma,\varphi)$-property} or a {\it $(\sigma,\varphi)$-group}. Moreover, everywhere below,
we denote by the symbols $N$ and $I$ the corresponding subgroups of a $(\sigma,\varphi)$-group
$G$ and denote the images $\sigma(H)$ and $\varphi(H)$ of the subgroups $H$ in $G$ by the
symbols $H^{\sigma}$ and $H^{\varphi}$.
\par
It is easily shown that for each coset $Ig$ there exists exactly one element of
$\varphi(G_{\sigma})$. Indeed, if $Ig=Ih$, then $g^{-1}g^{\sigma}=h^{-1}h^{\sigma}$. Hence,
$$
g=\varphi(g^{-1}g^{\sigma})=\varphi(h^{-1}h^{\sigma})=h.
$$
Denote the group $G$ regarded as a punctured symmetric space with the product $x.y=xy^{-1}x$
and the basic point 1 by the symbol $G^{+}$. Obviously, the restriction of the homomorphism
$\varphi$ of $N$ to $G_{\sigma}$ induces a homomorphism $G_{\sigma}\to G^{+}$ of the punctured
symmetric spaces. Let us an abstract symmetric space $S$ isomorphic to $\varphi(G_{\sigma})$,
and let $P:S\to\varphi(G_{\sigma})$ be the corresponding isomorphism. One can readily see
that, for every coset $Ig$, where $g\in N^{\varphi}$, there is exactly one element $x\in S$
such that
$$
Ig\cap\varphi(G_{\sigma})=\{P_x\}.
$$
(Here and everywhere below, for the sake of brevity, we use the symbol $P_x$ instead of
$P(x)$.) Hence, we can define a permutation representation of the group $N^{\varphi}$ on $S$
by setting $xg=z$ if $IP_{x}g=IP_{z}$ for $g\in N^{\varphi}$.
\medskip
{\bf Theorem 1.} {\sl If a group $G$ has the $(\sigma,\varphi)$-property, then the binary
composition $xy=xP_{y}$ equips $S$ with the structure of a Bol loop (of a Moufang loop
provided that $\varphi$ is an monomorphism).}
\medskip
{\bf Proof.} Let us fix a basic point $e$ in the space $S$. Since every homomorphism of
punctured spaces takes the basic point to the basic point, it follows that $P_{e}=1$ and
$P_{x^{-1}}=P_{x}^{-1}$. Therefore,
$$
ex=x=xe\qquad\text{and}\qquad (yx)x^{-1}=y.
$$
Further, $P_{x.y^{-1}}=P_{x}P_{y}P_{x}$. Therefore,
$$
x.y^{-1}=(xy)x,
$$
and hence
$$
z((xy)x)=((zx)y)x.
$$
We have thus proved that $S$ is a right Bol loop with the identity element. Let us now show
that the equation $ax=b$ has a unique solution in the loop $S$. Indeed, if a solution exists,
then it is unique, because
$$
a^{-1}((ax)a)=xa\qquad\text{and}\qquad x=(a^{-1}(ba))a^{-1}.
$$
On the other hand,
$$
a((a^{-1}(ba))a^{-1})=(ba)a^{-1}=b.
$$
Thus, a solution of the equation $ax=b$ exists. Hence $S$ is a Bol loop.
\par
Let $\varphi: N \to G$ is an isomorphism of the group $N$ into the group $G$. Then the mapping
$L:S\to G^+$ defined by the relation $L_{x}=P_{x}^{-1}P_{x}^{\sigma}$ is an isomorphism of the
punctured symmetric spaces. Therefore, $L_{x^{-1}}=L_{x}^{-1}$ and
$$
(P_{x}L_{y})^{-1}(P_{x}L_{y})^{\sigma}=P_{(y^{-1}x)y^{-1}}^{-1}P_{(y^{-1}x)y^{-1}}^{\sigma}.
$$
Hence,
$$
IP_{x}L_{y}=IP_{(y^{-1}x)y^{-1}}.
$$
Thus, the representation $x\to xL_{y}=(y^{-1}x)y^{-1}$ of the group $N$ on $S$ is defined. On
the other hand,
$$
P_{x^{-1}}^{\sigma}P_{x}=P_{x}P_{x^{-1}}^{\sigma}\qquad\text{and}\qquad
yP_{x}L_{x}=yL_{x}P_{x}.
$$
Therefore,
$$
(x^{-1}(yx))x^{-1}=x^{-1}y\qquad\text{and}\qquad x(x^{-1}y)=y.
$$
This implies that the Bol loop $S$ is a loop with inversion, and thus a Moufang loop. This
completes the proof of the theorem.
\smallskip
It is easy to prove that the generators $P_{x}$ and $R_{x}=P_{x^{-1}}^{\sigma}$ of any
$(\sigma,\varphi)$-group $G$ satisfy the following identities:
$$
\gather P_{e}=1,\quad R_{e}=1,
\tag 1\\
P_{x^{-1}}=P_{x}^{-1},\quad R_{x^{-1}}=R_{x}^{-1},
\tag 2\\
P_{(xy)x}=P_{x}P_{y}P_{x},\quad R_{(xy)x}=R_{x}R_{y}R_{x},
\tag 3\\
P_xP_yP_{xy}^{-1}=R_x^{-1}R_y^{-1}R_{xy}. \tag 4
\endgather
$$
Indeed, relations (1)--(3) immediately follow from Theorem 1. To prove relations (4), it
suffice to note that
$$
(P_xP_y)^{-1}(P_xP_y)^{\sigma}=P_{xy}^{-1}P_{xy}^{\sigma}.
$$
\medskip
{\bf Corollary 1.} {\sl If $G$ is a group with the $(\sigma,\varphi)$-property, then the
subgroup $I$ of $N^{\varphi}$ is generated by the elements $P_{x,y}=P_{x}P_{y}P_{xy}^{-1}$,
there $x,y\in S$. In particular, $I=\{g\in N^{\varphi}\mid eg=e\}$.}
\medskip
{\bf Proof.} Let $H=\langle P_{x,y}\mid x,y\in S\rangle$ and $K=\{g\in N^{\varphi}\mid
eg=e\}$, where $e$ is the identity element of $S$. Obviously, $H\subseteq I\subseteq K$. Let
us show that $K\subseteq H$. Indeed, $N^{\varphi}=\langle P_{x}\mid x\in S\rangle$. Therefore,
every element of $K$ can be represented as a word of the form $W=P_{x_{1}}P_{x_{2}}\dots
P_{x_{n}}$. If $n=1$, then $W=P_{e}\in H$. If $n>1$, then
$W=P_{x_{1},x_{2}}P_{x_{1}x_{2}}\dots P_{x_{n}}$, and the assertion is proved by an obvious
induction on the length of the word $W$.
\smallskip
The $(\sigma,\varphi)$-group $G$ and the Bol loop $S=S(G)$ constructed in Theorem 1 are said
to be {\it associated}. If $H$ is a subgroup of the group $G$, then we set
$$
H_{G}=\bigcap_{g\in G}gHg^{-1}.
$$
Obviously, $H_{G}$ is the largest normal subgroup of $G$ belonging to $H$. A subgroup $H$ is
said to be {\it $\varphi$-admissible} if
$$
(N\cap H)^{\varphi}\subseteq N^{\varphi}\cap H.
$$
A normal subgroup $K$ of the group $H$ is said to be {\it $\sigma$-admissible} if
$K^{\sigma}=K$ and the induced automorphism $\bar\sigma$ of $H/K$ is not identity. The
$\sigma$-admissible and $\varphi$-admissible normal subgroup $H$ of $G$ is said to be {\it
$(\sigma,\varphi)$-admissible}. If a group $G$ contains no proper
$(\sigma,\varphi)$-admissible normal subgroups, then this group is said to be {\it
$(\sigma,\varphi)$-simple}.
\medskip
{\bf Corollary 2.} {\sl Let $H$ be a $\sigma$-admissible normal subgroup of a
$(\sigma,\varphi)$-group $G$. The quotient group $\bar G=G/H$ has the
$(\bar\sigma,\bar\varphi)$-property if and only if $H$ is a $\varphi$-admissible subgroup of
$G$. The loop $\bar S$ associated with the group $\bar G$ is a homomorphic image of $S=S(G)$,
and $\bar S\simeq S$ if and only if $N^{\varphi}\cap H\subseteq I_{N^{\varphi}}$}.
\medskip
{\bf Proof.} Let $K=N\cap H$ and $M=N^{\varphi}\cap H$. Suppose $K^{\varphi}\subseteq M$.
Since $N/K\simeq N^{\varphi}/K^{\varphi}$, it follows that the homomorphism $G\to \bar G$
induces the natural homomorphism $\bar\varphi: N/K\to N^{\varphi}/M$ with the kernel
$M/K^{\varphi}$. It is obvious that $\bar G$ is a $(\bar\sigma,\bar\varphi)$-group.
Conversely, it follows from the existence of a homomorphism $\bar\varphi$ that its kernel
satisfies the relation $Ker(\bar\varphi)=M/K^{\varphi}$. Thus, $K^{\varphi}\subseteq M$. To
construct a homomorphism $S\to\bar S$, one can use the map $P_{x}\to P_{x}M$. In this case it
is obvious that $S\simeq\bar S$ if and only if the elements of $M$ trivially act on $S$, or
equivalently, on the cosets of $I$. This completes the proof of the assertion.
\smallskip
Let $S$ be a Bol loop associated with the $(\sigma,\varphi)$-group $G$, and let
$$
T_x=L_xR_x^{-1},\quad R_{x,y}=R_{x}R_{y}R_{xy}^{-1}, \quad L_{x,y}=L_{x}L_{y}L_{yx}^{-1},
$$
where $L_x=P_x^{-1}R_x^{-1}$ and $x,y\in S$. Using the obvious identity
$L_x^{\sigma}=L_x^{-1}$ it is easy to show that the set of all $L_x$ generates a
$\sigma$-invariant subgroup $N$ in $G$. Consider the homomorphic images $N^{\bar\varphi}$ and
$\bar I$ of $N$ and $N_L=\langle L_{x,y}\mid x,y\in S\rangle$ respectively arising as the map
$\bar\varphi=\sigma\varphi\sigma$. Since $L_x^{\sigma}=L_x^{-1}$, $\varphi(L_x)=P_x$ and
$P_x^{\sigma}=R_x^{-1}$, it follows that
$$
N^{\bar\varphi}=\langle R_x\mid x\in S\rangle,\quad \bar I=\langle R_{x,y}\mid x,y\in
S\rangle,
$$
and these groups are isomorphic to the groups $N^{\varphi}$ and $I$ respectively. Suppose
$\bar S$ is a symmetric space isomorphic to $S$. As above, we can  define a binary composition
on $\bar S$ setting $xy=z$ if $\bar IR_xR_y=\bar IR_{z}$. It is readily seen that this
composition equip $\bar S$ with the structure of a Bol loop. Obviously, the loops $\bar S$ and
$S$ are isomorphic. Finally, let
$$
J=\langle T_x, R_{x,y}, L_{x,y}\mid x,y\in S\rangle.
$$
\medskip
{\bf Proposition 1.} {\sl If $S$ is a Bol loop associated with the $(\sigma,\varphi)$-group
$G$, then the multiplication group $Gr(S)$ of $S$ is a homomorphic image of $G$. The kernel of
the homomorphism is $J_{G}$.}
\medskip
{\bf Proof.} It can easily be checked that
$$
JL_xR_y=JR_xR_y=JR_{xy}=JL_yL_x=JR_yL_x.
$$
In addition, $L_x^{-1}=R_xP_x$ and $P_x=L_{x^{-1}}R_{x^{-1}}$. Therefore,
$$
JL_x^{-1}R_y=JR_{x^{-1}y}\qquad\text{and}\qquad JR_xL_y^{-1}=JR_{(y^{-1}(xy))y^{-1}}.
$$
Since any element of $G$ is represented as a word of $R_x$ and $L_x$, we prove, by induction,
that every coset of $J$ has the element $R_x$ for some $x\in S$. Further, if $JR_x=JR_y$, then
$R_{xy^{-1}}\in J$. On the other hand, the intersection $J\cap \bar N^{\varphi}=\bar I$.
Therefore $R_{xy^{-1}}\in \bar I$, and hence $P_{xy^{-1}}\in I$. It is possible only if $x=y$.
Hence, every coset of $J$ in $G$ has exactly one element  $R_x$ for every $x\in S$.
\par
Thus, we can define a binary composition on the set of cosets of $J$ setting $xy=z$ if
$JR_xR_y=JR_{z}$ for all $x,y\in S$. Noting that $JR_x=JL_x$, we have
$$
JR_yL_x=JR_{xy}\quad\text{or}\quad yL_x=xy.
$$
Therefore the multiplication group $Gr(S)$ of $S$ is a homomorphic image of $G$. The kernel of
the homomorphism is the set of elements which induce the trivial permutation on $S$, or
equivalently, on cosets of $J$. Obviously, this kernel is the subgroup $J_{G}$. This completes
the proof of Proposition 1.
\smallskip
Suppose $J^{\sigma}\varSigma$ is a semidirect product of $J^{\sigma}$ and
$\varSigma=\langle\sigma\rangle$. Since the intersection $J^{\sigma}\cap N^{\varphi}=I$, it
follows that $(J^{\sigma}\varSigma)_{G}\cap N^{\varphi}= I_{N^{\varphi}}$. Therefore,
$(J^{\sigma}\varSigma)_{G}$ is the $(\sigma,\varphi)$-admissible normal subgroup of $G$. Using
Corollary 2 to Theorem 1, we have
\medskip
{\bf Corollary.} {\sl Let $H$ be a normal $(\sigma,\varphi)$-admissible subgroup of the
$(\sigma,\varphi)$-group $G$. Then the loops $S(G/H)$ and $S(G)$ are isomorphic if and only if
$H\subseteq (J^{\sigma}\varSigma)_{G}$}.
\bigskip

{\bf 3. Universality of the construction}
\medskip
Let $S$ be a Bol loop. Denote by $\widetilde G(S)$ the group presented by the generators
$P_{x}$ and $R_{x}$, where $x\in S$, and by the defining relations (1)--(4). Obviously,
$\widetilde G(S)$ admits an involutory automorphism $\sigma$ defined by the mapping
$P_{x}^{\sigma}=R_{x}^{-1}$ and $R_{x}^{\sigma}=P_{x}^{-1}$. One can present the group
$\widetilde G(S)$ in the form of a free product
$$
\widetilde G(S)=(A*B; A_{1}=B_{1},\psi)
$$
of the groups
$$
A=\langle P_{x}\mid x\in S\rangle\quad\text{and}\quad B=\langle R_{x}\mid x\in S\rangle
$$
with the subgroups
$$
A_{1}=\langle P_{x}P_{y}P_{xy}^{-1}\mid x,y\in S\rangle\quad\text{and} \quad B_{1}=\langle
R_{x}^{-1}R_{y}^{-1}R_{xy}\mid x,y\in S\rangle
$$
amalgamated according to the isomorphism $\psi: A\to B$. In this case, the automorphism
$\sigma$ of the group $\widetilde G(S)$ can be defined by setting $\sigma=\psi$ on $A$ and
$\sigma=\psi^{-1}$ on $B$.
\medskip
{\bf Theorem 2.} {\sl The group $\widetilde G(S)$ is a $(\sigma,\varphi)$-group, and
$S(\widetilde G(S))=S$. Any other $(\sigma,\varphi)$-group $G$ such that $S(G)=S$ is a
homomorphic image of the group $\widetilde G(S)$}.
\medskip
{\bf Proof.} Let $\widetilde N$ be the subgroup of $\widetilde G=\widetilde G(S)$ generated by
all elements of the form $L_{x}=P_{x}^{-1}R_{x}^{-1}$, where $x\in S$. Represent the relation
(4) in the form
$$
\align
R_{y}L_{x}R_{y}^{-1}&=L_{y}^{-1}L_{xy},\\
P_{y}^{-1}L_{x}P_{y}&=L_{xy}L_{y}^{-1},
\endalign
$$
and note that $L_{x}^{\sigma}=L_{x}^{-1}$. It is clear that $\widetilde N$ is a
$\sigma$-invariant normal subgroup of the group $\widetilde G$. Let $W$ be an element of
$\widetilde G$ represented as a word in the generators $P_{x}$ and $R_{x}$, where $x\in S$.
Then
$$
\align
(WP_x)^{-1}(WP_x)^{\sigma}&=P_x^{-1}W^{-1}W^{\sigma}P_xL_x,\\
(WR_x)^{-1}(WR_x)^{\sigma}&=R_x^{-1}W^{-1}W^{\sigma}R_xL^{-1}_{x^{-1}}.
\endalign
$$
Using induction on the length of the word $W$, we can prove that $\widetilde N=\langle
g^{-1}g^{\sigma}\mid g\in\widetilde G\rangle$. Since $A_1P_xP_y=A_1P_{xy}$, it is obvious that
every right coset of the subgroup $A_1$ in the group $A$ contains an element $P_x$ for some
$x\in S$.
\par
Now, let us consider the homomorphism $\pi$ of $\widetilde G$ onto $A$ such that $\pi(a)=a$
and $\pi(b)=b^{\sigma}$ for all $a\in A$ and $b\in B$. It is obvious that $Ker(\pi)=\widetilde
N$. Therefore the crossings $\widetilde N\cap A$ and $\widetilde N\cap B$ are trivial.
Hence~(see~[8]), $\widetilde N$ is a free group. On the other hand, there is a homomorphism of
the group $A$ onto the group of right multiplications of the loop $S$ such that every
generator $P_{x}$ is mapped to the right multiplication operator on $x$ in $S$. Therefore,
$P_x\in A_1$ only if $x=e$, there $e$ is the identity element of $S$, and $P_xP_y\in A_1$ only
if $y=x^{-1}$. Hence the subgroup $\widetilde N$ is freely generated by the elements $L_x$,
where $x\in S$ and $x\ne e$, and therefore the mapping $L_{x}\to P_{x}$ defines a homomorphism
$\varphi$ of $\widetilde N$ onto $A$. Thus, $\widetilde G$ is a $(\sigma,\varphi)$-group. It
then follows from Theorem 1 that $S(\widetilde G(S))=S$.
\par
Further, let $G$ be another $(\sigma,\varphi)$-group such that $S(G)=S$. Then $G$ has the set
of generators $P_{x}$ and $R_{x}$ ($x\in S$) satisfying the relations (1)--(4), and therefore
$G$ is a homomorphic image of $\widetilde G$. This completes the proof of the theorem.
\medskip
{\bf Corollary 1.} {\sl For any Bol loop $S$, there is an associated group $G_0=G_{0}(S)$
(which is unique up to isomorphism) such that any other group $G$ associated with $S$ has a
homomorphic image equal to $G_{0}$. Moreover, $(J^{\sigma}_0\varSigma)_{G_0}=1$.}
\medskip
{\bf Proof.} Let $\widetilde G =\widetilde G(S)$ and $G_{0}=\widetilde G/( \widetilde
J^{\sigma}\varSigma)_{\widetilde G}$. By Theorem 2, the groups $G$ and $G_{0}$ are homomorphic
images of $\widetilde G$. Suppose $G=\widetilde G/K$ for a $(\sigma,\varphi)$-admissible
normal subgroup $K$. By Corollary to Proposition 1, $K$ is contained in  $(\widetilde
J^{\sigma}\varSigma)_{\widetilde G}$. Since $(\widetilde J^{\sigma}\varSigma)_{\widetilde G}$
is a maximum $(\sigma,\varphi)$-admissible normal subgroup of $\widetilde G$, it follows that
$G_{0}$ is a minimal group associated with the loop $S$. Since $(\widetilde
J^{\sigma}\varSigma)_{\widetilde G}$ is a prototype of $(J^{\sigma}_0\varSigma)_{G_0}$ in
$\widetilde G$, it follows that   $(J^{\sigma}_0\varSigma)_{G_0}=1$. The corollary is proved.
\medskip
{\bf Corollary 2.} {\sl A Bol loop $S$ is simple if and only if the associated group
$G_{0}(S)$ is a $(\sigma,\varphi)$-simple. In addition, if the loop $S$ is strongly simple,
then the subgroup $N_{0}^{\varphi}$ of the group $G_{0}(S)$ is a simple group.}
\medskip
{\bf Proof.} By Corollary 1 to Theorem 2, it follows that $(J^{\sigma}_0\varSigma)_{G_0}=1$.
Therefore, by Corollary to Proposition 1, if $H$ is a proper $(\sigma,\varphi)$-admissible
normal subgroup of $G_{0}(S)$, then $G_{0}(S)/H$ is a group with the
$(\sigma,\varphi)$-property and its associated loop is a proper homomorphic image of $S$.
Conversely, let $S_{1}$ be a proper normal subloop of $S$ and $\phi:S\to S/S_1$ the induced
homomorphism. Then the maps $P_{x}\to P_{\phi(x)}$ and $R_{x}\to R_{\phi(x)}$ can be extended
to a homomorphism, provided that the image of an identity (1)--(4) in $\widetilde G(S)$ is an
identity in $\widetilde G(S/S_{1})$, and this readily verified. Therefore there exists the
homomorphism $\widetilde G(S)\to \widetilde G(S/S_{1})$. This homomorphism induces the
homomorphism $G_{0}(S)\to G_{0}(S/S_{1})$ with a $(\sigma,\varphi)$-admissible kernel.
Therefore the group $G_{0}$ is not $(\sigma,\varphi)$-simple.
\par
Further, suppose $S$ is a strongly simple Bol loop and $K$ is a normal subgroup of the group
$N_{0}^{\varphi}$. Then the Proposition 1 required the condition $K\subseteq (J_0)_{G}\cap
I_0$. On the other hand, the group $I_0$ is a set of fixed points of $\sigma$. Hence,
$K^{\sigma}\subseteq (J^{\sigma}_0\varSigma)_{G_0}$. By Corollary 1 to Theorem 2, this implies
that $K=1$. The corollary is proved.
\smallskip
Denote by $Gr_{r}(S)$ and $Gr_{l}(S)$ the right and left multiplication groups of the Bol loop
$S$ respectively. The following assertion describes a structure of the multiplication group
$Gr(S)$ in case if the loop $S$ is strongly simple.
\medskip
{\bf Proposition 2.} {\sl Let $S$ be a strongly simple Bol loop, and $H$ is a normal subgroup
of $Gr(S)$. Then $Gr(S)$ is a $\varphi$-admissible group, and one of the following assertions
are holds: \roster
\item"{$i)$}"
$Gr(S)$ is a simple group and $S$ is a Moufang loop;
\item"{$ii)$}"
$H=Gr_{l}(S)$, and $Gr(S)=H\leftthreetimes Gr_{r}(S)$ is a semidirect product of the
subgroups;
\item"{$iii)$}"
$H\simeq Gr_{r}(S)$, and $Gr(S)=Gr_{l}(S)\times H$ is a direct product of the subgroups.
\endroster}
{\bf Proof.} If $G_0=G_0(S)$ is a simple group, then $N_0=G_0$, and hence the subgroup
$Ker(\varphi)$ of $N_0$ is a normal subgroup of the group $G_0$. Obviously, $Ker(\varphi)\ne
G_0$. Therefore, $Ker(\varphi)=1$. Using Theorem 1, we prove that $S$ is a Moufang loop.
\par
Let $H$ be a nonidentity proper normal subgroup of $G_{0}$, and let $K=H\cap N_{0}$. If
$K\nsubseteq Ker(\varphi)$, then $K^{\varphi}$ is a nonidentity normal subgroup of
$N_{0}^{\varphi}$. By Corollary 2 to Theorem 2, it is possible only if $K=N_{0}$. Therefore,
$N_{0}\subseteq H$. On the other hand,
$$
G_0/H\simeq (G_0/N_0)/(H/N_0)\simeq N_{0}^{\varphi}/(H/N_0).
$$
Since the group $N_{0}^{\varphi}$ is simple, it follows that $H=N_{0}$. In addition, it is
obvious that $N_0\cap N_0^{\varphi}=1$. On the other hand, it is clear that the group $G_0$ is
generated by the subgroups $N_0$ and $N_0^{\varphi}$. Therefore $G_0$ is a semidirect product
of the subgroups.
\par
Let $K\subseteq Ker(\varphi)$. Since $Ker(\varphi)_{G_0}$ is the largest normal subgroup of
$G_0$ belonging to $Ker(\varphi)$, it follows that $K\subseteq Ker(\varphi)_{G_0}$. On the
other hand, $HN_0$ is a normal subgroup of $G_0$. As was shown above, it is not a proper
normal subgroup of $G_{0}$. Hence, $HN_0=G_0$.
\par
Thus, the quotient group $G_{0}/Ker(\varphi)_{G_0}$ is a simple or semisimple group. On the
other hand, Corollary 1 to Theorem 1 claims that $I_{0}\subseteq J_{0}^{\varphi}$. Therefore
$Ker(\varphi)\subseteq J_{0}$, and hence $Ker(\varphi)_{G_{0}}\subseteq (J_{0})_{G_{0}}$.
Using Proposition 1, we prove the assertion.
\bigskip

{\bf 4. Finite strongly simple Bol loops}
\medskip

Recall (see [8]) that a group $G$ is said to be {\it residually finite} if the intersection of
all its normal subgroups of finite index is the identity group. A group $G$ is residually
finite if and only if, for every element $g\in G$, $g\ne 1$ there is a homomorphism of $G$
into a finite group that takes $g$ to a nonidentity element. As is known (see~[9]), in
particular, every free product of two finite groups with an amalgamated subgroup is a
residually finite group.
\medskip
{\bf Theorem 3.} {\sl If $S$ is a finite Bol loop of order $n$, then $\widetilde G=\widetilde
G(S)$ is a residually finite group. In particular, the group $G_0$ is finite.}
\medskip
{\bf Proof.} It is obvious that the group $\widetilde N^{\varphi}$ admits an involutory
automorphism $\theta$ such that $P_{x}^{\theta}=P_{x}^{-1}$ for all $x\in S$. Consider the
semidirect products $\widetilde N^{\varphi}\langle\theta\rangle$. By the relations (1)--(3),
it follows that
$$
(P_{x}\theta)^{\theta}=P_{x^{-1}}\theta,\qquad
(P_{x}\theta)^{P_{y}}=P_{(y^{-1}x)y^{-1}}\theta.
$$
These relations imply that the group $\widetilde N^{\varphi}\langle\theta\rangle$ is generated
by a finite set of elements each of which is of finite order and has finitely many conjugate
elements.
\par
Let $|S|=n$. Then $X=\{P_{x_{1}}\theta,\dots, P_{x_{n}}\theta\}$ is the set of generators of
the group $\widetilde N^{\varphi}\langle\theta\rangle$ and of the elements conjugate to the
generators. Let us show that every element of $\widetilde N^{\varphi}\langle\theta\rangle$ can
be represented by a word $W$ of length $\leq n$. Indeed, $(P_{x_{i}}\theta)^{2}=1$. Therefore,
$$
W=(\dots P_{x_{i}}\theta W_{1}P_{x_{i}}\theta\dots) =(\dots
(P_{x_{i}}\theta)^{-1}W_{1}P_{x_{i}}\theta\dots).
$$
Since, when conjugating $P_{x_{j}}\theta$ by $P_{x_{i}}\theta$, we obtain one of the elements
in the list $X$ again, it follows that the length of the word $W$ can be reduced. After
finitely steps we obtain a representation of $W$ in the form of a word of length $\leq n$. It
is easily shown that the order
$$
|\widetilde N^{\varphi}|\leq\frac12\sum^{n}_{k=1}\frac{n!}{(n-k)!}<\frac12en!,
$$
where $e$ is the base of the natural logarithms. Since every free product of two finite groups
with an amalgamated subgroup is a residually finite group, it follows that the first part of
the assertion is proved.
\par
Thus, $N_0$ is a finite-index subgroup of $G_0$, and $Ker(\varphi)$ is a finite-index subgroup
of $N_0$. This implies that $K=Ker(\varphi)_{G_0}$ is a finite-index normal subgroup of $G_0$.
Hence, $G_0$ is a finite group. This completes the proof of the theorem.
\smallskip
By analogy with the notion of an algebraic group we define an {\it algebraic Bol loop}. One is
a Bol loop $S$ equipped with an algebraic variety structure such that the product $S\times
S\to S$ and the inversion $S\to S$ are regular mappings (morphisms) of the algebraic
varieties.
\medskip
{\bf Theorem 4.} {\sl If $S$ is an connected algebraic strongly simple Bol loop, and its
multiplication group $Gr(S)$ is an algebraic group, then $S$ is a Moufang loop. In particular,
every  connected finite strongly simple Bol loop is Moufang.}
\medskip
{\bf Proof.} As is well known~(see e.g.~[10,11]), there exists two basic type of connected
algebraic groups: {\it Abelian varieties} and {\it linear (or affine) algebraic groups}. Any
Abelian variety is an Abelian group. Any linear algebraic group are a subgroup of the group
$GL(n,k)$. These two classes of  groups have the trivial intersection. It is known that:
\roster
\item"{$i)$}"
any connected algebraic group $G$ has unique normal linear algebraic subgroup $H$ such that
the quotient group $G/H$ is an Abelian variety;
\item"{$ii)$}"
any linear algebraic group $G$ has a maximum connected solvable normal subgroup ({\it a
radical}) $K$ such that the quotient group $G/K$ is semisimple.
\endroster
\par
Let $S$ be an connected algebraic strongly simple Bol loop, and $Gr(S)$ be an algebraic group.
Obviously, $Gr_{r}(S)$ and $Gr_{r}(S)^{\sigma}$ are connected subgroups of $Gr(S)$. Let $G(M)$
be an intersection of all closed subgroups of $Gr(S)$ containing the set $M=Gr_{r}(S)\cup
Gr_{r}(S)^{\sigma}$. Then~(see~[10]) $G(M)$ is a connected algebraic subgroup of $Gr(S)$.
Since the set $M$ generates the group $Gr(S)$, it follows that $Gr(S)$ is a connected
algebraic group.
\par
Suppose $H$ is a normal linear algebraic subgroup of $Gr(S)$ such that $Gr(S)/H$ is an Abelian
manifold. It follows from Propositions 2 that either $H=1$ or $H=Gr_{l}(S)$ or $H\simeq
Gr_{r}(S)$ or $H=Gr(S)$. If $H=1$, then $Gr(S)$ is an Abelian group, and hence $S$ is a simple
Abelian group. If $H=Gr_{l}(S)$ or $H\simeq Gr_{r}(S)$, then $Gr_{r}(S)$ is a linear algebraic
group and an Abelian variety simultaneously. It is possible only if $S$ is the identity group.
\par
Let $Gr(S)$ be a linear algebraic group, and let $K$ be a radical of $Gr(S)$. If $K\ne 1$,
then either $K=Gr(S)$ or $K=Gr_{l}(S)$ or $K\simeq Gr_{r}(S)$. In any case $Gr_{r}(S)$ is a
solvable group. Since $Gr_{r}(S)$ is a simple group, it follows that it is an Abelian group,
and hence $S$ is a simple commutative Moufang loop. If $K=1$, then $Gr(S)$ is a simple group,
and $S$ is a Moufang loop.
\par
Finally, let $S$ be an finite connected algebraic strongly simple Bol loop. It follows from
Theorem 3 that the group $Gr(S)$ is finite, and hence one may be embedded into the linear
algebraic group $GL(n,k)$. As above, this implies that $Gr(S)$ is a connected finite group,
and hence $S$ is a finite simple Moufang loop. This completes the proof of Theorem 4.
\medskip
{\bf Corollary.} {\sl If the right multiplication group $Gr_{r}(S)$ of a Bol loop $S$ is a
finite simple group of Lie type, then $S$ is a finite simple Moufang loop.}
\medskip
{\bf Proof.} Obviously, $S$ is a finite loop. Therefore $Gr(S)$ is a finite group, and hence
one may be embedded into $GL(n,k)$. As in the previous theorem, this implies that $Gr(S)$ is a
connected algebraic group, and hence $S$ is a finite simple Moufang loop.
\medskip
{\it Remark 1.} Apparently, Theorem 4 may be extended to the class of all connected finite
simple Bol loop. The key to the proof is the assertion that the subgroup $N_{0}^{\varphi}$ for
such loop is a simple group. However, we have not any proof of the assertion.
\medskip
{\it Remark 2.} It is known~[12] that the group $Aut(S)$ of all invertible morphisms of a loop
$S$ into itself has a natural structure of an algebraic group. Obviously, $Gr(S)$ is a
subgroup of $Aut(S)$. However, the group $Aut(S)$ can have infinite number of connected
components, and hence the assertions i) and ii) for $Aut(S)$ may not be true. The example that
was constructed in~[6] shows that there exists a strongly simple non-Moufang Bol loop with a
non-algebraic multiplication group.
\bigskip

\centerline {\bf Acknowledgements}
\medskip
The author is thankful to professor D.I.~Moldavanskii for discussion of the content of the
paper and to the referee for many useful remarks that enabled me to improve the exposition.
The research supported by RFBR Grant 06-02-16140.
\bigskip

\centerline {\bf References}
\medskip

\noindent\item{1.} Robinson,~D.A. Bol loops. Trans\. Amer\. Math\. Soc. {\bf 1966}, {\it 123},
341-354.

\noindent\item{2.} Gorenstein,~D.; Lyons,~R.; Solomon,~R. {\it The classification of finite
simple groups}. Amer\. Math\. Soc.; Surveys and Monographs, 1995; Vol.~40.

\noindent\item{3.} Doro,~S. Simple Moufang loops. Math\. Proc\. Camb\. Phil\. Soc. {\bf 1978},
{\it 83}, 377-392.

\noindent\item{4.} Liebeck,~M.W. The classification of finite simple Moufang loops. Math\.
Proc\. Camb\. Phil\. Soc. {\bf 1987}, {\it 102}, 33-47.

\noindent\item{5.} Bruck,~R. {\it A Survey of Binary Systems}; Springer-Verlag: Berlin, 1971.

\noindent\item{6.} Kiechle,~H.; Kinyon,~M. Infinite simple Bol loops. Comment\. Math\. Univ\.
Caroline. {\bf 2004}, {\it 45}, 275-278.

\noindent\item{7.} Loos,~O. {\it Symmetric Spaces}; New York-Amsterdam, 1969.

\noindent\item{8.} Magnus,~W,; Karrass,~A.; Solitar,~D. {\it Combinatorial Group Theory}; New
York-London-Sydney, 1966.

\noindent\item{9.} Baumslag,~G. On the residual finiteness of generalized free products of
nilpotent groups, Trans\. Amer\. Math\. Soc. {\bf 1963}, {\it 106}, 193-209.

\noindent\item{10.} Borel,~A. {\it Linear Algebraic Groups}; New York: Benjamin W.A., 1969.

\noindent\item{11.} Rosenlicht,~M. Some basic theorems on algebraic groups. Amer\. J\. Math.
{\bf 1956},  {\it 78}, 401-443.

\noindent\item{12.} Matsusaka~T. Polarized varieties fields of moduli and generalized Kummer
varieties of polarized abelian varieties. Amer\. J\. Math. {\bf 1958}, {\it 80}, 45-82.

\end